\documentclass [11pt,oneside,a4paper,mathscr]{amsart}

\usepackage{amscd,amsmath,amssymb,euscript}

\newtheorem{thm}{Theorem}[section]
\newtheorem{cor}[thm]{Corollary}
\newtheorem{lemma}[thm]{Lemma}
\newtheorem{prop}[thm]{Proposition}
\newtheorem{defn}[thm]{Definition}

\theoremstyle{definition}
\newtheorem{example}[thm]{Example}
\newenvironment{pf}{\paragraph{Proof}}{\par\medskip}

\newcommand{\mat}[4]{\left( \begin{array}{cc} #1 & #2 \\ #3 & #4
\end{array} \right)}
\newcommand{\col}[2]{\left( \begin{array}{c} #1 \\ #2 \end{array} \right)}
\renewcommand{\S}{\mathcal{S}}
\newcommand{\id}{\operatorname{id}}
\newcommand{\la}{\lambda}
\newcommand{\K}{K}
\newcommand{\SL}{\operatorname{SL}}

\newcommand{\curlym}{\phi}
\newcommand{\cP}{\P}
\newcommand{\I}{\mathcal I}
\newcommand{\NN}{\mathcal N}
\newcommand{\spell}{{\ell}}
\newcommand{\PP}{\mathbb{P}}
\newcommand{\ch}{\operatorname{ch}}
\newcommand{\End}{\operatorname{End}}
\newcommand{\Hilb}{\operatorname{Hilb}}
\newcommand{\ellsur}{\pi\colon X\to S}

\newcommand{\C}{\mathbb C}
\newcommand{\Z}{{\mathbb Z}}
\newcommand{\ha}{\frac{1}{2}}
\newcommand{\Oh}{\mathcal O}

\newcommand{\De}{\Delta}
\newcommand{\D}{\operatorname{D}}
\newcommand{\Ext}{\operatorname{Ext}}
\newcommand{\blob}{{\scriptscriptstyle\bullet}}
\newcommand{\eu}{\operatorname{\chi}}
\newcommand{\Hom}{\operatorname{Hom}}
\newcommand{\lHom}{\operatorname{\mathcal H\mathit{om}}}
\newcommand{\lra}{\longrightarrow}
\newcommand{\embed}{\hookrightarrow}
\newcommand{\lRa}[1]{\stackrel{#1}{\longrightarrow}}
\newcommand{\lla}{\longleftarrow}
\newcommand{\lLa}[1]{\stackrel{\,#1}{\lla}}
\newcommand{\isom}{\cong} 
\renewcommand{\L}{\mathbf L}
\newcommand{\R}{\mathbf R}
\newcommand{\dual}{\vee}
\newcommand{\tensor}{\otimes}
\newcommand{\Ltensor}{\stackrel{\mathbf L}{\tensor}}
\newcommand{\OO}{\Oh}
\newcommand{\m}{\mathfrak{m}}
\newcommand{\n}{\mathfrak{n}}
\newcommand{\Tor}{\operatorname{Tor}}
\newcommand{\depth}{\operatorname{depth}}
\newcommand{\hd}{\operatorname{hd}}
\newcommand{\supp}{\operatorname{supp}}
\newcommand{\codim}{\operatorname{codim}}
\newcommand{\im}{\operatorname{im}}
\renewcommand{\P}{\mathcal P}
\newcommand{\Q}{\mathcal Q}
\newcommand{\E}{\mathcal E}
\newcommand{\hpi}{\Hat{\pi}}
\newcommand{\M}{\operatorname{\mathcal M}}
\newcommand{\Spec}{\operatorname{Spec}}

\newcommand{\cl}{\operatorname{c}}

\newcommand{\rk}{\operatorname{r}}
\newcommand{\Pl}{\operatorname{P}}
\newcommand{\QQ}{\mathbb {Q}}
\newcommand{\TT}{\operatorname{T}}

\title{Fourier-Mukai transforms for K3 and elliptic fibrations}
\author{Tom Bridgeland \and Antony Maciocia}

\begin{document}

\begin{abstract}
Given a non-singular variety with a K3 fibration $\pi\colon X\to S$
we construct dual fibrations $\hpi\colon Y\to S$ by replacing
each fibre $X_s$ of $\pi$ by a two-dimensional moduli space of stable
sheaves on $X_s$.
In certain cases we prove that the resulting scheme $Y$ is a non-singular
variety and construct an equivalence of derived categories
of coherent sheaves $\Phi\colon \D(Y)\to\D(X)$. Our methods also apply to elliptic and
abelian surface fibrations. As an application we use the
equivalences $\Phi$ to relate moduli spaces of stable bundles
on elliptic threefolds to Hilbert schemes of curves.
\end{abstract}

\maketitle

\section{Introduction}

Moduli spaces of stable vector bundles on projective surfaces
have been the subject of a great deal of research in recent years,
and much is now known about them. In contrast, moduli spaces of
bundles on higher-dimensional varieties are still rather poorly
understood.
In this paper we extend the theory of Fourier-Mukai transforms,
a useful tool in the study of moduli spaces of bundles on surfaces,
to an interesting class of higher-dimensional varieties,
namely those with K3 or elliptic fibrations. 

\subsection{}

Let $\pi\colon X\to S$ be a  Calabi-Yau fibration,
that is a connected morphism of non-singular projective varieties
whose general fibre has trivial canonical bundle. Suppose
that $\pi$ has relative dimension one or two, so that the generic
fibre is an abelian surface, a K3 surface or an elliptic curve. We
shall also fix some polarisation $\ell$ of $X$.

Let $Y$ be a component of the relative moduli space of stable
sheaves on the fibration $\pi$.
Points of $Y$ represent stable sheaves
supported on the fibres of $\pi$,
and there is a  natural map $\hpi\colon Y\to S$ sending
a sheaf supported on the fibre
$\pi^{-1}(s)$ to the corresponding point $s\in S$.

These relative moduli spaces have been studied by other authors,
notably R. Thomas \cite{Th}. We shall denote them by $\M^{\ell}(X/S)$.
Here we shall be mainly interested in cases when
$Y$ has the same dimension as $X$,
and is fine, in that there is a universal sheaf on $Y\times X$.
We shall show that $Y$ is then a non-singular
projective variety and the morphism $\hpi\colon Y\to S$
is another Calabi-Yau fibration, which we refer to 
as a Mukai dual of the fibration $\pi$.

Clearly, one would expect the geometry of the two spaces $X$ and $Y$
to be closely related, and indeed, this
turns out to be the case.
Thus, the Hodge numbers of the two spaces
often coincide, and at a deeper level,
one finds that certain
moduli spaces of sheaves on $X$ and $Y$
can be identified.
One of the  aims of this paper is to study this geometrical relationship
in more detail.
It turns out that a good way to do this
is to observe that the universal sheaf on $Y\times X$
induces a natural equivalence
between the derived categories of coherent sheaves on the two spaces.
Such equivalences are called Fourier-Mukai (FM) transforms,
the first example being due to Mukai \cite{Muk1}.

The construction outlined above was applied to
surfaces fibred by elliptic curves in \cite{Br1}.
The resulting FM transforms 
were used to show that the general component of the
moduli space of stable sheaves on an elliptic surface is
birational to a Hilbert schemes of points.
The same transforms
were also  used to describe a type of mirror symmetry for string
theories compactified on
elliptically fibred K3 surfaces \cite{BB}.

In this paper we extend these techniques
to higher dimensions by generalising the results of \cite{Br1} to
varieties fibred by K3 or abelian surfaces, or by elliptic curves,
over bases of arbitrary dimension. 
We give some applications of the resulting FM transforms
to moduli spaces of
stable bundles on threefolds. One might also expect
these transforms to play a r{\^o}le in the
mathematical description of mirror symmetry
for hyperk{\"a}hler fourfolds. 

Note that Calabi-Yau fibrations
can have singular fibres, which in general
are both non-reduced and reducible. Describing the possible
forms for such fibres is a difficult problem.
An interesting aspect of our methods
is that we do not need to know this information. We prove a kind
of `removable singularities' result, which, put roughly, states that
providing the singular fibres have sufficiently high codimension
it is enough to understand the case of non-singular fibres.

\subsection{}
We shall only consider fibrations
which satisfy the following
minimality condition.

\begin{defn}
\label{def}
In this paper a Calabi-Yau fibration is a morphism of non-singular projective
varieties $\pi\colon X\to S$ whose general fibre is a variety with
trivial canonical bundle and such that $\K_X\cdot C=0$ for any curve
$C\subset X$ contained in a
fibre of $\pi$.
\end{defn}

Note that such fibrations need not be flat, in fact, a morphism
of non-singular projective varieties is flat precisely
when it is equidimensional.
Note also that flat Calabi-Yau fibrations of relative dimension at
most two (with our definition)
are of one of three types: elliptic fibrations, K3 fibrations or
abelian surface fibrations. Our main result is 

\begin{thm}
\label{main}
Let $\pi\colon X\to S$ be a flat Calabi-Yau fibration
of relative dimension at most two. Take a polarization
$\ell$ on $X$ and
let $Y$ be an irreducible component of $\M^{\ell}(X/S)$.
Suppose that
$Y$ is fine and of the same dimension as $X$ and that the morphism $\hpi\colon Y\to S$
is equidimensional.
Let $\P$ be a
universal sheaf on $Y\times X$. Then $Y$ is non-singular,
$\hpi$ is a flat Calabi-Yau fibration of the same type as $\pi$,
and the functor
\[\Phi^{\P}_{Y\to X}(-)=\R\pi_{X,*}(\P\tensor\pi_Y^*(-)),\]
where $Y\lLa{\pi_Y}Y\times X\lRa{\pi_X} X$ are the projection maps,
is an equivalence of derived categories $\D(Y)\to \D(X)$.
\end{thm}

There are well-known methods for finding fine components of $\M^{\ell}(X/S)$
of the right dimension which we describe in Section \ref{sixone} below.
However, the condition that $\hpi$ is equidimensional
is somewhat unsatisfactory, since it is difficult
to check in practice. The problem is that there
conceivably exist Calabi-Yau fibrations with
certain fibres $X_s$ whose singularities are so severe that the
corresponding moduli space $Y_s$ has
dimension larger than that of $X_s$.
Nonetheless the authors do not know of an example when
this problem actually occurs,
and in particular, we shall show that
when $X$ has dimension at most three,
and the morphism $\hpi$ is surjective,
the condition that $\hpi$ is equidimensional
is automatic.

The main difficulty in the proof of Theorem \ref{main} is in showing
that the
dual variety $Y$ is non-singular. This point did not arise in \cite{Br1},
because the dimension of the
tangent space to $Y$ at any point
could be calculated directly using the Riemann-Roch formula.
As we remarked above, much less is known about the geometry of moduli
spaces of sheaves on higher-dimensional varieties,
so the statement that $Y$ is non-singular is surprisingly strong.

The proof we give uses a deep result
in commutative algebra known as the intersection theorem.
Since our argument might seem rather complicated at first sight,
it may be
worth noting that the methods we develop here
have also been used to
show that three-dimensional Gorenstein quotient singularities have crepant
resolutions \cite{BKR},
a fact only known previously through case-by-case analysis.

\subsection{}
The statement that two varieties $X$ and $Y$ have equivalent
derived categories has some strong geometrical consequences. Firstly,
certain topological invariants of $X$ and $Y$ should coincide.
One expects, for example, that the Hodge numbers of the two spaces are
equal.
This statement should follow from Barannikov and Kontsevich's work on
$A^{\infty}$ deformations of derived categories \cite{BK}.
In the meantime we have the following result, first proved in special
cases by R. Thomas \cite{Th}.

\begin{prop}
\label{wow}
In the situation of Theorem \ref{main}, if
$X$ is a Calabi-Yau threefold,
then so is $Y$, and then
$h^{p,q}(X)=h^{p,q}(Y)$ for all $p,q$.
\end{prop}

A more important consequence of an equivalence of derived categories
is that various moduli spaces of sheaves on the two varieties become
identified. Thus we have

\begin{thm}
\label{moduli}
Let $X$ be a threefold with a flat elliptic fibration
$\pi\colon X\to S$ and let $\hpi\colon Y\to S$ be a Mukai dual fibration
as in Theorem \ref{main}. Let $\NN$ be a connected component of
the moduli space of rank one, torsion-free sheaves on $Y$.
Then there is a polarisation $\ell$ of $X$ and a connected component
$\M$ of the moduli space of stable torsion-free sheaves on $X$
with respect to $\ell$
which is isomorphic to $\NN$.
\end{thm}

This is a pretty result, but in practice one would like
to be able to explicitly describe a given moduli
space on $X$.
This is a more delicate problem which we leave for future research.
Here we content ourselves with an example; we show that a particular
moduli space of stable bundles on a
non-singular anticanonical divisor in $\PP^2\times\PP^2$ is isomorphic
to $\PP^2$.

\subsection*{Plan of the paper}
We start in Sections 2 and 3 with some preliminaries on moduli spaces of stable
sheaves and Fourier-Mukai transforms respectively.
The intersection theorem is stated in Section 4 and interpreted geometrically
in Section 5. These ideas are used in Section 6 to
prove the `removable singularities' result, Theorem \ref{biggy}, which,
in turn, is used to prove
Theorem \ref{main} in Section 7. Sections 8 and 9
contain a more detailed treatment of the case of elliptically fibred
threefolds
and we conclude in Section 10 with the example of a non-singular anticanonical
divisor in $\PP^2\times\PP^2$.

\subsection*{Acknowledgements}
This paper has benefited from our discussions with many
mathematicians. In particular, it is a pleasure to thank
Richard Thomas who read an earlier draft and suggested
several improvements, and Mark Gross
who helped with the geometry of elliptic threefolds.
Thanks also to an anonymous reviewer who pointed
out an error in an earlier version. Finally we should acknowledge our debt to
Shigeru Mukai, whose brilliant paper \cite{Muk3} contains not only the main results
we need concerning moduli spaces of sheaves on K3 surfaces, but also
the key ideas which allow one to construct Fourier-Mukai 
transforms. This research was carried out with the support of the
Engineering and Physical Sciences Research Council of Great Britain.

\subsection*{Conventions}
We work throughout in the category of schemes over $\C$. Points of a scheme
are always closed points. By a sheaf on a scheme $X$ we mean a coherent
$\OO_X$-module. The bounded derived category
of coherent sheaves on a scheme $X$ is denoted $\D(X)$. The $i$th
cohomology (resp. homology) sheaf of an object $E$ of $\D(X)$
is denoted $H^i(E)$ (resp. $H_i(E)$),
and an object of $\D(X)$ satisfying $H^i(E)=0$ for all $i\neq 0$
will be referred to simply as a sheaf.


\section{Moduli spaces of stable sheaves}

The concept of stability of sheaves was introduced by Mumford in the case of
bundles on curves and later extended to torsion-free sheaves on surfaces by
Gieseker and Maruyama.
More recently Simpson defined a general notion of stability for
sheaves on arbitrary projective schemes. In this section we summarise
some of Simpson's results. For proofs see \cite{Si}.

\subsection{}
Let $X$ be a projective scheme over $\C$ and let $\ell$ be a polarisation
of $X$. Thus $\ell$ is the first Chern class of an ample line bundle $L$ on $X$. 
The Hilbert polynomial of a sheaf $E$ on $X$
with respect to the polarisation $\ell$ is the polynomial function
$\Pl_E(n)=\eu(E\tensor L^{\tensor n})$.
The coefficients of this polynomial
are given in terms of the Chern classes of $E$ by the
Riemann-Roch theorem.
In particular the degree of $\Pl_E$ is the dimension
of the support of $E$. The unique rational multiple of $\Pl_E$ which
is monic is called the normalised Hilbert polynomial of $E$ and will be
denoted $\cP_E$. Polynomials of the same degree are ordered so that
$\cP_1<\cP_2$ precisely when $\cP_1(n)<\cP_2(n)$ for $n\gg 0$.

A sheaf $E$ on $X$ is said to be of pure dimension
if the support of any non-zero subsheaf of $E$ has the same
dimension as the support of $E$. For example, if $X$ is a variety and $E$
has positive rank this is just the condition that $E$ is torsion-free.
A sheaf $E$ is stable with respect to the polarisation $\ell$
if $E$ has pure dimension, and if for all proper subsheaves $A\subset E$
one has $\cP_A<\cP_E$. The usual properties of stability
carry over to this more general situation. For example, if $A$ and $B$ are
stable with $\cP_A\geq\cP_B$ then any non-zero map $A\to B$ is an isomorphism.
In particular, stable sheaves are simple.

There is a moduli space $\M^{\ell}(X)$
of stable sheaves on $X$, which in general has infinitely
many components, each of which is a quasi-projective scheme
which may be compactified by adding equivalence classes of semistable
sheaves.
An irreducible component $Y\subset\M^{\ell}(X)$ is said to be \emph{fine}
if $Y$ is projective (so that there is no semistable boundary),
and there exists a universal family of stable sheaves
$\{\P_y: y\in Y\}$ on $X$. In particular there is a
sheaf $\P$ on
$Y\times X$, flat over $Y$, such that for each $y\in Y$, $\P_y$ is the
stable sheaf on $X$ represented by the point $y$. Note that by the
universal
property, the family  $\{\P_y: y\in Y\}$ is complete,
that is for each point $y\in Y$ the Kodaira-Spencer map
$\TT_y Y\to \Ext^1_X(\P_y,\P_y)$
is an isomorphism.

Any given irreducible component $Y\subset\M^{\ell}(X)$ parameterises sheaves $E$
with fixed numerical invariants, and hence fixed Hilbert polynomial $\Pl$.
A result of Mukai \cite[Theorem A.6]{Muk3} implies that $Y$ is fine whenever
the integers
$\eu(E\tensor L^{\tensor n})$ have no common factor.

We now turn to relative moduli spaces.
Let $\pi\colon X\to S$ be a morphism of projective schemes,
and fix a polarization $\ell$ of $X$.
For each point $s\in S$, $\ell$ induces a polarization $\ell_s$ of the fibre
$X_s$ of $\pi$, defined by restricting an ample line bundle $L$
representing
$\ell$ to the subscheme $X_s$.
In this situation there is a relative moduli scheme
\[\hpi\,\colon \M^{\ell}(X/S)\lra S\]
whose fibre over a point $s\in S$ is naturally the moduli scheme
$\M^{\ell_s}(X_s)$ of stable sheaves on the fibre $X_s$,
with respect to the polarisation $\ell_s$.
Points of $\M^{\ell}(X/S)$ correspond to stable sheaves $E$ on $X$
whose scheme-theoretic support is contained in
some fibre $X_s$ of $\pi$. The morphism $\hpi$ takes the point
representing such a sheaf $E$ to the corresponding point $s\in S$. The
existence of such moduli schemes is due originally to Simpson (see
\cite{Si} or \cite[Section 4]{Th}) but is now standard.

As before, if $Y$ is a connected component of $\M^{\ell}(X/S)$
parameterising sheaves $E$ for which the integers
$\eu(E\tensor L^{\tensor n})$ have no common factor,
then $Y$ is fine, in that $Y$ is projective
and there is a universal sheaf $\P$ on $Y\times X$.
This universal sheaf is actually supported on the closed subscheme
$j\colon Y\times_S X\embed Y\times X$.
Thus there is a sheaf $\P^+$ on $Y\times_S X$, flat over $Y$,  such that
$\P=j_*\P^+$. For each point $y\in Y$ with $\hpi(y)=s$, 
one has $\P_y={i_y}_*\P^+_y$, where $\P^+_y$ is the restriction of
$\P^+$ to the subscheme
$i_y\colon\{y\}\times X_s\embed X$.

\subsection{}
We finish this section with two simple results about stable sheaves.
Both are presumably well-known, but we include
proofs since we were unable to find a suitable reference.

\begin{lemma}
\label{fib}
Let $\pi\colon X\to S$ be a morphism of non-singular projective
varieties of relative dimension one. Let $\ell$ be a polarisation of $X$ and let $\curlym$ 
be the pull-back of a polarisation from $S$.
Let $E$ be a torsion-free sheaf on $X$ whose restriction to the
general fibre of $\pi$ is stable. Then there exists
an integer $m_0$ such that for all $m\geq m_0$ the sheaf $E$ is
stable with respect to the polarisation $\ell+m\curlym$.  
\end{lemma}

\begin{pf}
Let $d$ be the dimension of $S$. For each proper subsheaf $A\subset E$ put
\[\De(A)=r(A)\cl_1(E)-r(E)\cl_1(A)\in H^2(X,\Z).\]
The sheaf $E$
is called $\mu$-stable with respect to a polarisation $\omega$ if
for all such $A$ one has $\De(A)\cdot\omega^d >0$. In particular
$\mu$-stable sheaves are stable. The lemma follows from the following
two facts. Firstly, for fixed $m$, the numbers
$\De(A)\cdot(\ell+m\curlym)^d$ are
bounded below. Secondly, the fact that the restriction of $E$ to the general
fibre of $\pi$ is stable implies that $\De(A)\cdot\curlym^d>0$.
\qed
\end{pf}

\begin{lemma}
\label{hilb}
Let $X$ be a non-singular, simply-connected projective variety.
Let $\NN$ denote the union of those
connected components of $\Hilb(X)$ which parametrise closed subschemes
of codimension greater than 1. Let $\M$ be the moduli space of (stable)
torsion-free sheaves with rank one and zero first Chern class. Then
there is a morphism of schemes $\alpha\colon\NN\to\M$ inducing a
bijection on closed points, which sends a subscheme $C\subset X$ to
the ideal sheaf $\I_C$. 
\end{lemma}

\begin{pf}
There is a universal sheaf $\P$ on $\NN\times X$, flat over $\NN$,
such that for each point $s\in\NN$, $\P_s$ is the ideal sheaf of the
corresponding subscheme $C\subset X$. The existence of the map
$\alpha$ follows from the universal property of $\M$.

A point of $\M$ corresponds to a torsion-free sheaf $E$ of rank 1 and zero
first Chern class. There is a short exact sequence
\[0\lra E\lra E^{**}\lra T\lra 0\]
where $T$ is supported in codimension $>1$ and
the double dual sheaf $E^{**}$ is reflexive of rank 1,
hence invertible \cite[Lemma II.1.1.15]{OSS}.
Since $X$ is simply-connected and $\cl_1(E^{**})=0$,
one has $E^{**}=\OO_X$,
so $E=\I_C$ for some subscheme $C\in\NN$.
Thus $\alpha$ is bijective on points.
\qed
\end{pf}

The authors do not know of any examples where the map $\alpha$ of the Lemma
fails to be an isomorphism.


\section{Fourier-Mukai transforms}

Let $X$ be a non-singular projective variety. The derived category
of $X$ is a triangulated category whose objects are complexes of
sheaves on $X$ with bounded and coherent homology sheaves. 
This category contains a large amount of geometrical information about $X$.
Indeed, if $\pm K_X$ is ample, Bondal and Orlov \cite{BO} 
show that $\D(X)$ considered up
to triangle-preserving equivalence
completely determines $X$. In general however, there exist
pairs of non-singular projective varieties $(X,Y)$ for which there are
triangle-preserving equivalences
$\Phi\colon\D(Y)\to\D(X).$
Such equivalences are called
Fourier-Mukai transforms.
In this section we shall describe
some of the geometrical consequences of a pair of varieties being
related in this way. First however, we address the problem of
constructing examples of FM transforms.

\subsection{Constructing  FM transforms}
A theorem of Orlov \cite{O} states that for any FM transform
$\Phi\colon\D(Y)\to\D(X),$
there is an object $\P$ of $\D(Y\times X)$ such that $\Phi$ is isomorphic
to the functor $\Phi^{\P}_{Y\to X}$ given by the
formula
\[\Phi^{\P}_{Y\to X}(-)=\R\pi_{X,*}(\P\Ltensor\pi_Y^*(-)),\]
where $Y\lLa{\pi_Y}Y\times X\lRa{\pi_X} X$ are the projection maps.
It remains to determine which objects $\P$ give equivalences.
Suppose for simplicity that $\P$ is a sheaf on $Y\times X$, flat over
$Y$ and put $\Phi=\Phi^{\P}_{Y\to X}$.

If $A$ and $B$ are objects of $\D(X)$ we put
\[\Hom_{\D(X)}^i(A,B)=\Hom_{\D(X)}(A,B[i])\]
where $[i]$ is the functor which shifts all complexes of $\D(X)$ to the
left by $i$ places. If $A$ and $B$ are sheaves on $X$, that is complexes whose
homology is concentrated in degree zero, one has
$\Hom^i_{\D(X)}(A,B)=\Ext^i_X(A,B).$

The key observation is that if $\Phi$ is a FM transform then
\[\Hom^i_{\D(X)}(\Phi(A),\Phi(B))=\Hom^i_{\D(Y)}(A,B).\]
This suggests the idea of thinking of $\Hom^{\blob}(-,-)$ as analogous to an
inner product, so that the functor $\Phi$ becomes an isometry.
Pursuing the analogy, note that for distinct points
$y_1,y_2$ of $Y$, one has
$\Ext^i_Y(\OO_{y_1},\OO_{y_2})=0$  for all $i$.
In this way (see also Lemma \ref{numpty}) the set of sheaves
$\{\OO_y:y\in Y\}$
can be thought of as an orthogonal
basis for $\D(Y)$. Applying $\Phi$
must give an orthogonal basis of sheaves $\{\P_y:y\in Y\}$ on $X$,
parameterised by $Y$. The point of the theorem below is that conversely,
such families gives rise to FM transforms.

\begin{thm}
\label{borlov}
Let $X$ and $Y$ be non-singular projective varieties of the same dimension,
and $\P$ a sheaf on $Y\times X$, flat over $Y$. Then the functor
\[\Phi^{\P}_{Y\to X}\colon \D(Y)\lra\D(X)\]
is an equivalence of categories if and only if, for any point $y\in Y$,
\[\End_X(\P_y)=\C\mbox{ and }\P_y\tensor\omega_X=\P_y,\]
and for any pair of distinct points $y_1$, $y_2$ of $Y$, and any
integer $i$,
\[\Ext^i_X(\P_{y_1},\P_{y_2})=0.\]

\end{thm}

\begin{pf}
A complete proof, based on ideas of Bondal, Orlov and Mukai appears in
\cite{Br2}. There is also a more general version, which we shall not need in
this paper, in which $\P$ is allowed
to be an arbitrary object of $\D(Y\times X)$.
\qed
\end{pf}

\subsection{Applications}

\label{adjy}

Suppose $\Phi=\Phi^{\P}_{Y\to X}\colon\D(Y)\to\D(X)$ is a FM transform.
The inverse of $\Phi$ is the functor
$\Psi=\Phi^{\Q}_{X\to Y}\colon\D(X)\to\D(Y)$
where $\Q$ is the object 
\[\Q=\R\lHom_{\OO_{Y\times X}}(\P,\pi_X^* \omega_X)[n],\]
and $n$ is the dimension of $X$ and $Y$.
Thus there are isomorphisms
\[\Psi\circ\Phi\isom\id_{\D(Y)},\qquad\Phi\circ\Psi\isom\id_{\D(X)}.\]
Using the Chern classes of the objects $\P$ and $\Q$ 
one can define a pair of mutually inverse
correspondences relating the cohomology of $X$ and $Y$.
This leads to the following result. For a proof see 
\cite[Theorem 4.23]{Th} or \cite[Theorem 4.9]{Muk2}.

\begin{prop}
\label{co}
Let $X$ and $Y$ be non-singular projective varieties. A FM transform
$\Phi\colon\D(Y)\to\D(X)$
induces isomorphisms
\[ \bigoplus_i H^{i,i+k}(Y) \lra \bigoplus_i H^{i,i+k}(X)\]
for each integer $k$.
\qed
\end{prop}

Proposition \ref{wow} follows easily from Proposition \ref{co} once one
observes that Serre duality requires that if
$\Phi\colon\D(Y)\to\D(X)$ is a FM transform
then $Y$ has trivial canonical
bundle precisely when $X$ does \cite[Lemma 2.1]{Br3}.

\medskip

The second application of FM transforms we wish to consider is to
moduli spaces of sheaves. Before explaining the general
method we should introduce some notation.
Let $\Phi\colon\D(Y)\to\D(X)$ be a FM transform. Given a sheaf
$E$ on $Y$ we write $\Phi^i(E)$ for the $i$th cohomology sheaf of the
object $\Phi(E)$ of $\D(X)$.

\begin{defn}
A sheaf $E$ on $Y$ is
said to be $\Phi$-WIT$_i$ if $\Phi^j(E)=0$ unless $j=i$. A sheaf $E$ is
$\Phi$-WIT if it is $\Phi$-WIT$_i$ for some $i$. One then refers
to the sheaf $\Hat{E}=\Phi^i(E)$ as the transform of $E$.
\end{defn}

The reason that FM transforms can be used to compute moduli
spaces of sheaves is that
they preserve families \cite[Proposition 4.2]{Br2}.
Thus if $\{\E_s:s\in S\}$ is a family of $\Phi$-WIT
sheaves on $Y$, then $\{\Hat{\E_s}:s\in S\}$ is a family
of sheaves on $X$. As a consequence one has the well-known

\begin{prop}
\label{jo}
Let $\{\E_s:s\in S\}$
be a complete family of pairwise non-isomorphic sheaves on $Y$,
over a connected, projective scheme $S$.
Suppose that for each point $s\in S$ the sheaf $\E_s$ is $\Phi$-WIT
and the transform $\Hat{\E_s}$ is stable
with respect to some polarisation $\ell$ of $X$. Then
the natural map $S\to\M^{\ell}(X)$
which sends a point $s\in S$ to the point
representing the sheaf $\Hat{\E_s}$ is
an isomorphism onto a
connected component of $\M^{\ell}(X)$. 
\end{prop}

\begin{pf}
The map exists by the universal property of
$\M^{\ell}(X)$ and the fact that $\{\Hat{\E_s}:s\in S\}$
is a family of stable
sheaves on $X$. Since this family is complete
this map induces isomorphisms on tangent spaces.
It is injective on points because the sheaves
$\Hat{\E_s}$ are pairwise non-isomorphic.
The result follows.
\qed
\end{pf}

In practice
one often takes the sheaves $\E_s$
to be torsion-free of rank $1$, since moduli spaces of such sheaves
can be related to Hilbert schemes via Lemma \ref{hilb}. Many moduli
spaces have been computed in this way.


\section{The Intersection Theorem}

Our proof of Theorem \ref{main} uses a celebrated result
of commutative algebra called the intersection theorem. The aim of this
section is to state and prove a strong version of this theorem
using Hochster's results on the existence of maximal Cohen-Macaulay modules.
Our main reference is \cite{Syz}.

We start by defining the depth of an arbitrary (not necessarily
finitely generated) module over a local ring.

\begin{defn}
Let $(A,\m)$ be a Noetherian local ring. The depth of an $A$-module $M$ is
defined to be
\[\depth_A(M)=\inf\{i\in\Z:\Ext^i_A(A/\m,M)\neq 0\}.\]
\end{defn}

In general modules may have infinite depth, but for finitely generated
modules the above definition is always finite and agrees with the
usual definition of depth in
terms of $A$-sequences. See \cite[Theorem 1.6]{Syz}.
 
A fairly immediate consequence of this definition is the following
acyclicity lemma of Peskine and Szpiro. For a proof see \cite[Lemma 1.3]{Syz}.

\begin{lemma}
Let $(A,\m)$ be a Noetherian local ring and
\[0\lra M_s\lra M_{s-1}\lra\cdots\lra M_0\lra 0,\]
a finite complex of $A$-modules such that each non-zero homology module
$H_i(M_{\blob})$ has depth $0$. Then
\[\depth_A(M_i)\geq i\qquad 0\leq i\leq s\]
implies that $H_i(M_{\blob})=0$ for all $i>0$, and
\[\depth_A(M_i)> i\qquad 0\leq i\leq s\]
implies that $M_{\blob}$ is exact, that is $H_i(M_{\blob})=0$ for all
$i$.
\qed
\end{lemma}

We can now state the main result of this section. The first part is
called the intersection theorem; we
reproduce the proof for the reader's convenience. The second
part is a slight strengthening which we prove using similar
methods. See \cite[Theorem 1.13]{Syz}.

\begin{thm}
\label{com}
Let $(A,\m)$ be a Noetherian local $\C$-algebra of dimension
$d$ and let
\[0\lra M_s\lra M_{s-1}\lra\cdots\lra M_0\lra 0,\]
be a non-exact complex of finitely-generated free $A$-modules such
that each homology module $H_i(M_{\blob})$ is a finite $A$-module.  Then
$s\geq d$. Furthermore, if $s=d$ and $H_0(M_{\blob})\isom A/\m$, then
\[H_i(M_{\blob})=0\mbox{ for all }i\neq 0,\]
and $A$ is regular.
\end{thm}

\begin{pf}
We may assume that $(A,\m)$ is complete.
Fixing a system of parameters $(x_1,\cdots ,x_d)$ for
$A$, let $R$ denote the complete regular local ring
\[R=\C[[x_1,\cdots x_d]],\]
and let $\n$ denote its maximal ideal. Note that $A$ is finite as an $R$-module.

Work of Griffith and Hochster \cite[Theorem 1.8]{Syz} shows that
in this situation there is a (not-necessarily
finitely-generated) $A$-module $C$ such that $C$ is free as an
$R$-module. It follows immediately that $\depth_A(C)=d$. One says that
$C$ is a balanced big Cohen-Macaulay module for $A$.

Let us put
\[N=A\tensor_R (R/\n).\]
Then $N$ is a finite $A$-module satisfying\footnote{Note added in 2019. In fact this  step is not clear, and there is therefore a significant gap in the published version of this paper. Fortunately Theorem \ref{com} is still correct; for a proof and further discussion see T. Bridgeland and S. Iyengar, `A criterion for regularity of local rings',  C. R. Math. Acad. Sci. Paris 342 (2006), no. 10, 723--726. The authors would like to thank Srikanth Iyengar for pointing out this mistake.}
\begin{equation}
\label{one}
\Tor_p^A(N,C)=0\mbox{ for all } p>0,
\end{equation}
and $N\tensor C\neq 0$. It follows from this that for any non-zero finite
length $A$-module $P$, the module $P \tensor
C$ is non-zero.

Suppose now that $s<d$, and tensor the complex $M_{\blob}$ by
$C$. Each term of the resulting complex $M_{\blob}\tensor C$ has depth $d$
because it is a direct
sum of finitely many copies of $C$. Applying the acyclicity lemma then
shows that $M_{\blob}\tensor C$ is exact. Consider the third-quadrant spectral sequence
\[E^{p,q}_2=\Tor_{-p}^A(H_{-q}(M_{\blob}),C)\implies
H_{-(p+q)}(M_{\blob}\tensor C).\]
If the complex $M_{\blob}$ is not exact, let $q_0$ be the least
integer such that $H_q(M_{\blob})$ is non-zero. Then
the term $E^{0,-q_0}_2$ of the sequence is non-zero and survives to infinity,
which is a contradiction. This proves the first part.

Assume now that $s=d$. By the acyclicity lemma again, the complex
$M_{\blob}\tensor C$ is exact except at
the right-hand end, so the spectral sequence $E^{p,q}_2$ converges to
0 unless $p+q=0$. In particular we must have $E^{-1,0}_2=0$, so since
we have assumed $H_0(M_{\blob})=A/\m$, we have
\[\Tor_1^A(A/\m,C)=0.\]
Using the fact that the finite $A$-module $N$ has a composition
series whose factors are isomorphic to $A/\m$, and (\ref{one}) above,
it is easy to see that for any finite $A$-module $P$
\[\Tor_p^A(P,C)=0\mbox{ for all }p>0.\]
Thus the spectral sequence degenerates and we can conclude that
$H_i(M_{\blob})=0$ for $i>0$, and so $M_{\blob}$ is a finite free resolution of
$A/\m$. It follows that $A$ is regular.
\qed
\end{pf}


\section{Support and homological dimension}

In this section we translate the intersection theorem into geometrical
language whence it becomes a statement relating the homological dimension of
a complex of sheaves to the codimension of the support of its
homology sheaves. Throughout $X$ denotes an arbitrary scheme of
finite type over $\C$.

\begin{defn}
The \emph{support} of an object $E$ of $\D(X)$, written $\supp(E)$,
is the union of the supports of the homology
sheaves $H_i(E)$ of $E$. It is a closed subset of $X$.
\end{defn}

\begin{defn}
Given a non-zero object $E$ of $\D(X)$, the  \emph{homological dimension} of
$E$, written $\hd(E)$, is equal
to the smallest integer $s$ such that $E$ is quasi-isomorphic to a complex
of locally free $\OO_X$-modules of length $s$. If no such integer
exists we put $\hd(E)=\infty$.
\end{defn}

The following two simple results allow one to calculate both the
support and homological dimension of a given object $E$ of $\D(X)$
from a knowledge of the vector spaces
\begin{equation}
\label{tb}
H_i(E\Ltensor\OO_x)=\Hom^i_{\D(X)}(E,\OO_x)^{\dual},
\end{equation}
as $x$ ranges over the points of $X$. To check that the above two spaces
really are equal, let $f:\{x\}\embed X$ be
the inclusion map, and apply \cite[Proposition II.5.6, Corollary II.5.11]{Ha}.

\begin{lemma}
\label{numpty}
Let $E$ be an object of $\D(X)$, and fix a point $x\in X$. Then
\[x\in\supp(E)\iff\exists i\in\Z:\Hom_{\D(X)}^i(E,\OO_x)\neq 0.\]
\end{lemma}

\begin{pf}
There is a spectral sequence
\[E^{p,q}_2=\Ext^p_X(H_q(E),\OO_x)\implies\Hom_{\D(X)}^{p+q}(E,\OO_x).\]
If $x$ lies in the support of $E$, let $q_0$ be the minimal value of
$q$ such that $x$ is contained in the support of the homology sheaf
$H_q(E)$.
Then there
is a non-zero element of $E^{0,q_0}_2$
which survives to give a non-zero element of
$\Hom_{\D(X)}^{q_0}(E,\OO_x)$. The converse is clear.
\qed
\end{pf}

\begin{prop}
\label{tordim}
Let $X$ be a quasi-projective scheme, take a non-zero object $E$ of $\D(X)$
and let $s\geq 0$ be an integer
such that for all points $x\in X$,
\[H_i(E\Ltensor\OO_x)=0\mbox{ unless }0\leq i\leq s.\]
Then $E$ is quasi-isomorphic to a complex of locally free
sheaves of the form
\[0\lra L_s\lra L_{s-1}\lra \cdots\lra L_0\lra 0.\]
In particular, $E$ has homological dimension at most $s$.
\end{prop}

\begin{pf}
If $s=0$
the result follows from \cite[Lemma 4.3]{Br1}, so assume $s$
positive. Using (\ref{tb}) and applying the argument used to prove the
last lemma, one sees that $H_i(E)=0$ for all $i<0$. Since $X$ is
quasi-projective, every coherent sheaf on $X$ is the quotient of a
finite rank, locally free sheaf, so the dual of
\cite[Lemma I.4.6]{Ha} implies that $E$ is quasi-isomorphic to a
complex of locally free sheaves on $X$ of the form
\[
\begin{CD}
\cdots @>>> L_i @>{d_i}>> L_{i-1} @>{d_{i-1}}>> L_{i-2} @>>> \cdots @>d_1>>
 L_0 @>>> 0.
\end{CD}\]
Consider the object $F$ of $\D(X)$ defined by the
truncated complex
\[\begin{CD}
\cdots @>>> L_{s+3} @>d_{s+3}>> L_{s+2} @>d_{s+2}>> L_{s+1} @>d_{s+1}
>> L_s @>>> 0.
\end{CD}\]
Applying the functor $-\tensor\OO_x$ and comparing with $E$, one sees that for
any $x\in X$, $H_i(F\Ltensor \OO_x)=0$
for $i>s$ and $i<s$. Applying the $s=0$ case, this implies that
$F$ has homological dimension 0. It follows that $H_i(E)=0$ for all
$i>s$, and that $L_s/\im(d_{s+1})$ is locally free.

Consider the short exact sequence of complexes
\[0\lra A\lra L\lra B\lra 0,\]
where $A$ is the complex
\[\begin{CD}
\cdots @>>> L_{s+2}@>d_{s+2}>>L_{s+1}@>>>\im (d_{s+1})@>>> 0,\end{CD}\]
and $B$ is the length $s$ complex of locally free sheaves
\[\begin{CD}0@>>> L_s/\im (d_{s+1})@>>>
L_{s-1}@>d_{s-1}>>\cdots@>>>L_0@>>> 0.\end{CD}\]
It is enough to show
that $A$ is quasi-isomorphic to zero. But this is clear since
$H_i(A)=H_i(L)=0$ for $i>s$.
\qed
\end{pf}

In geometrical language the intersection theorem becomes

\begin{cor}
\label{vanish}
Let $X$ be a scheme of finite type over $\C$ and $E$ a non-trivial object of
$\D(X)$. Then for any irreducible component $\Gamma$ of $\supp(E)$ one
has an inequality
\[\codim(\Gamma)\leq \hd(E).\]
\end{cor}

\begin{pf}
We may assume that $E$ is a complex of locally free sheaves of finite
length $s\geq 0$. Let $E_0$ be the
restriction of $E$ to the
affine subscheme $\Spec(A)$ of $X$, where $A$ is the local
$\C$-algebra $\OO_{X,{\Gamma}}$.
Then $E_0$ is non-trivial, and each homology sheaf $H_i(E_0)$ is a finite
$A$-module \cite[Corollary 2.18]{Ei}. Furthermore, the dimension of
$A$ is equal to the codimension of $\Gamma$. Theorem \ref{com}
now gives the result.
\qed
\end{pf}

The second part of Theorem \ref{com}, together with Lemma \ref{numpty}
and Proposition \ref{tordim} gives

\begin{cor}
\label{smooth}
Let $X$ be an irreducible quasi-projective scheme of dimension $n$ over $\C$
and fix a point $x\in
X$. Suppose that there is an object
$E$ of
$\D(X)$ such that for any point $z\in X$, and any integer $i$,
\[\Hom_{\D(X)}^i(E,\OO_z)=0 \mbox{ unless }z=x \mbox{ and }0\leq i\leq n.\]
Suppose also that $H_0(E)\isom\OO_{x}$. Then $X$ is non-singular at $x$,
and $E\isom\OO_{x}$.
\qed
\end{cor}


\section{A removable singularities result}

In this section we use the results of the last two sections
to prove the following strengthening
of Theorem \ref{borlov}. This is the result which allows us
to deal with singular fibres of Calabi-Yau fibrations.

\begin{thm}
\label{biggy}
Let $X$ be a non-singular projective variety of dimension $n$ and let
$\{\P_y:y\in Y\}$ be a complete family of simple sheaves on $X$ parameterised
by an irreducible projective scheme $Y$ of dimension $n$. 
Suppose that
\[\Hom_X(\P_{y_1},\P_{y_2})=0\]
for any distinct points $y_1,y_2\in Y$ and that the closed subscheme
\[\Gamma(\P)=\{(y_1,y_2)\in Y\times Y: \Ext^i_X(\P_{y_1},\P_{y_2})\neq
0\mbox{ for some }i\in\Z\}\]
of $Y\times Y$ has dimension at most $n+1$. Suppose also that
$\P_y\tensor\omega_X\isom\P_y$ for all $y\in Y$.
Then $Y$ is a non-singular
variety and the functor
\[\Phi^{\P}_{Y\to X}\colon \D(Y)\lra \D(X)\]
is an equivalence of categories.
\end{thm}

\begin{pf}
The sheaf $\P$ is flat over $Y$, and $X$ is
non-singular, so given a point $(y,x)\in Y\times X$, the complex
\[\P\Ltensor\OO_{(y,x)}\isom \P_y\Ltensor\OO_x\]
has bounded homology. It follows from Proposition \ref{tordim} that $\P$ has
finite homological dimension, and in particular, the object
\[\P^{\dual}=\R\lHom_{\OO_{Y\times X}}(\P,\OO_{Y\times X})\]
has bounded homology. By Grothendieck-Verdier duality the functor
$\Phi=\Phi^{\P}_{Y\to X}$ has a left adjoint
\[\Psi\colon\D(X)\lra\D(Y).\]
By composition of correspondences \cite[Proposition 1.3]{Muk1}, there is an object $\S$ of
$\D(Y\times Y)$ such that there is an isomorphism of
functors
\[\Psi\circ \Phi(-)\isom \R\pi_{2,*}(\S\Ltensor\pi_1^*(-)),\]
where $\pi_1$ and $\pi_2$ are the projections of the product $Y\times Y$ onto its
two factors.

For any point $y\in Y$, the derived restriction of the object $\S$ to
the subscheme $\{y\}\times Y$ is just $\Psi\Phi\OO_y$. It follows that for any pair
of points $(y_1,y_2)$ in $Y\times Y$,
\begin{equation}
\label{two}
H_p(\S\Ltensor\OO_{(y_1, y_2)})^{\dual}=\Hom^p_{\D(Y)}(\Psi\Phi\OO_{y_1},\OO_{y_2})=
\Ext^p_X(\P_{y_1},\P_{y_2})
\end{equation}
using the adjunction $\Psi\vdash \Phi$.
Since $X$ is non-singular these groups vanish unless $0\leq p\leq
n$.

If $y_1\neq y_2$ are distinct points of $Y$,
\[\Ext^n_X(\P_{y_2},\P_{y_1})=\Hom_X(\P_{y_1},\P_{y_2})=0,\]
so if we define $\E$ to be the restriction of $\S$ to the complement
of the diagonal in $Y\times
Y$, Proposition \ref{tordim} implies that $\E$ has
homological dimension $n-2$. By hypothesis the support of $\E$ has
codimension at least $n-1$, so applying Proposition \ref{vanish} shows that $\E\isom
0$. It follows that the support of $\S$ is the diagonal, and hence the groups
(\ref{two}) vanish unless $y_1=y_2$.

Fix a point $y\in Y$. For any other point $z\in Y$,
\[\Hom_{\D(Y)}^i(\Psi\Phi\OO_y,\OO_z)=\Ext^i_X(\P_{y},\P_z),\]
and these groups vanish unless $y=z$ and $0\leq i\leq n$. We claim that
$$H_0(\Psi\Phi\OO_y)=\OO_y.$$
Assuming this for the moment, note that
Corollary \ref{smooth} now implies that $Y$ is non-singular, and so Theorem
\ref{biggy} follows from Theorem \ref{borlov} (or from a simple piece
of category theory \cite[Theorem 2.3]{Br2}).

To prove the claim note first that there is a unique map
$\Psi\Phi\OO_y\to \OO_y$, so we obtain a triangle
\begin{equation}
\label{rt}
C\lra \Psi\Phi\OO_y\lra \OO_y\lra C[1]
\end{equation}
for some object $C$ of $\D(Y)$, which must be supported at $y$.
Applying the functor
$\Hom_{\D(Y)}(-,\OO_y)$ and using the adjunction
$\Psi\vdash \Phi$ gives a long exact sequence
\[0\lra\Hom_{\D(Y)}^0(\OO_y,\OO_y)\lra
\Hom_{\D(X)}^0(\Phi\OO_y,\Phi\OO_y)\lra \Hom_{\D(Y)}^0(C,\OO_y)\]
\[\lra \Hom_{\D(Y)}^1(\OO_y,\OO_y)\lRa{\epsilon}
\Hom_{\D(X)}^1(\Phi\OO_y,\Phi\OO_y)\lra \cdots \]
The homomorphism $\epsilon$ is just the Kodaira-Spencer map for the
family $\{\P_y:y\in Y\}$, \cite[Lemma 4.4]{Br2}, which is an isomorphism by
assumption. Since $\Phi\OO_y$ is simple, it follows that
\[\Hom_{\D(Y)}^i(C,\OO_y)=0\mbox{ for all }i\leq 0,\]
and so, by the argument of Lemma \ref{numpty}, $H_i(C)=0$ for all
 $i\leq 0$.
Taking
homology of the triangle (\ref{rt}) shows that $H_0(\Psi\Phi\OO_y)=\OO_y$
as claimed.
\qed
\end{pf}

Note that the condition $\Hom_X(\P_{y_1},\P_{y_2})=0$ for all distinct
$y_1$ and $y_2$ is automatic for \textit{stable} moduli spaces which
we will be considering in what follows.


\section{FM transforms for Calabi-Yau fibrations}

In this section we use Theorem \ref{biggy} to
prove our main result, Theorem \ref{main}. First
we explain how to construct components
of the moduli space of stable sheaves
on a Calabi-Yau fibration which satisfy the hypotheses of the
theorem. The reader should compare the treatment given by Thomas
\cite[Section 4]{Th}.

\label{mainbit}

\subsection{}
\label{sixone}

Let $\pi\colon X\to S$ be a Calabi-Yau fibration,
fix a polarisation $\ell$ of $X$
and let $X_s$ be a non-singular fibre of $\pi$ 
with its induced polarisation $\ell_s$.
Let $E$ be a sheaf on $X_s$, stable with respect to
$\ell_s$, and let $Y$ be the
irreducible component of $\M^{\ell}(X/S)$ containing a point representing
the sheaf $E$.

Suppose first that $\pi$ has relative dimension two.
By the Riemann-Roch theorem on $X_s$, the
Hilbert polynomial of $E$ with respect to $\ell$ is
\[\Pl_E(n)=\eu(E)+\cl_1(E)\cdot n\spell_s+\ha \rk(E)(n\spell_s)^2.\]
If the integers $\Pl_E(n)$ have no common factor
then the moduli
space $Y$ is fine.

There is an open subset
$U\subset S$, such that all the fibres $X_s$ over points $s\in U$ are
non-singular. The fibre $Y_s$ of $\hpi\colon Y\to S$ over a point
$s\in U$ is a component of
the moduli space of stable sheaves on the surface $X_s$.
Results of Mukai \cite[Theorems 0.1, 1.17]{Muk2} imply that
$Y_s$ is non-singular and either empty or of dimension
\[\cl_1(E)^2-2\rk(E)\eu(E)+\eu(\OO_{X_s})\rk(E)^2+2.\]
Thus the moduli space $Y$ has the same dimension as $X$ providing
$\hpi$ is surjective and
\[2\rk(E)\eu(E)=\eu(\OO_{X_s})\rk(E)^2+\cl_1(E)^2.\]
In this case the general fibre of $\hpi$ will be a Calabi-Yau surface
of the same type as the general fibre of $\pi$.

The condition that $\hpi$ is surjective is satisfied providing
$\cl_1(E)$ is the restriction of the first Chern class of some
line bundle on $X$. Note that in the case when $X$ and $Y$ have dimension
three, this is enough to show that $\hpi$ is equidimensional, because
the assumption that $Y$ is irreducible prevents the dimension
of the fibres of $\hpi$ from jumping.

\begin{example}
Suppose $\pi\colon X\to S$ is a K3 fibration, take a non-singular fibre
$X_s$ of $\pi$ and let $E$ be the ideal sheaf of a point on
$X_s$. Then $\eu(E)=1$ and for any polarisation $\ell$
of $X$ the component $Y\subset\M^{\ell}(X/S)$ containing $E$ is fine
and has the same dimension as $X$. The resulting FM transforms
are relative versions of Mukai's reflection functor
\cite[Section 2]{Muk3}.
\end{example}

Suppose now that $\pi$ has relative dimension
one. As before, let $E$ be a stable sheaf on a
non-singular fibre $X_s$ of $\pi$, necessarily an elliptic curve.
The numerical invariants of $E$ are
simply its rank $r$ and degree $d$. The condition that the irreducible
component
$Y$ of $\M^{\ell}(X/S)$ containing  $E$ be fine is just that $d$ is
coprime to both $r$ and the degree of the induced polarization
$\ell_s$. Whenever this happens the resulting moduli space has the
same dimension as $X$ and is elliptically fibred over $S$.
We shall give examples later.

\subsection{}

The main difficulty in the proof of Theorem \ref{main} is
dealing with singularities in the fibres of $\pi$, 
so we first consider the case when
all the fibres of $\pi$ are non-singular. To be definite
we shall assume that $\pi$ has relative dimension two.

Take a polarization $\ell$ of $X$ and let $Y$ be a fine component
of the relative moduli scheme $\M^{\ell}(X/S)$ of the same dimension
as $X$. Assume that $\hpi\colon Y\to S$ is flat. Each fibre of
$\hpi$ is non-singular, so $Y$ is a non-singular projective
variety. Let $\P^+$ be a universal sheaf on $Y\times_S X$.
Extending by zero we obtain a sheaf  $\P$ on $Y\times X$.
By definition $\P$ is flat over $Y$.

Fix a point $s\in S$, and let $i\colon X_s\embed X$ be the inclusion
morphism. For any point $y\in Y$ with $\hpi(y)=s$, one has $\P_y=i_*
\P^+_y$. Thus to check the first condition of Theorem \ref{borlov}
it is enough to show that $i^*(\omega_X)=\OO_{X_s}$.
This is immediate from the adjunction formula,
since $X_s$ has trivial canonical bundle and trivial normal bundle.

For the second condition we may assume that the two sheaves
$\P_{y_1}$, $\P_{y_2}$ are supported on the same fibre of $\ellsur$,
otherwise their $\Ext$-groups trivially vanish. Then
\[\Ext_X^i(\P_{y_1},\P_{y_1})=\Ext^i_X(i_*\P^+_{y_1}, i_*\P^+_{y_2})
=\Hom^i_{X_s}(\L i^* i_*\P^+_{y_1},\P^+_{y_2}),\]
by the adjunction $\L i^*\vdash i_*$. Furthermore
\[\L i^* i_*(\P^+_{y_1})=\P^+_{y_1}\Ltensor \L i^* i_*(\OO_{X_s})\]
by the projection formula. Since $S$ is non-singular, one can write down
the Koszul resolution for $\OO_s$ on $S$ and pull-back via $\pi$
to obtain a locally free resolution of $\OO_{X_s}$ on $X$. This gives
\[\L_q i^* i_*\P^+_{y_1}=\P^+_{y_1}\tensor\bigwedge^q\OO_{X_s}^{\oplus {m}},\]
for $0\leq q\leq s$, where $m=\dim S$ is the dimension of the base.
Now there is a spectral sequence
\[E^{p,q}_2=\Ext^p_{X_s}(\L_q i^* i_*\P^+_{y_1},\P^+_{y_2})\implies
\Ext^{p+q}_X(\P_{y_1},\P_{y_2}),\]
so it is enough to know that
\[\Ext^i_{X_s}(\P^+_{y_1}, \P^+_{y_2})=0,\]
for all $i$. Since the set $\{\P^+_y\colon y\in Y_s\}$ is a two-dimensional
moduli space of stable sheaves on the surface $X_s$, this
is an immediate consequence of the Riemann-Roch theorem
\cite[Proposition 3.12]{Muk3}.

\subsection{}

\label{bl}
We now prove Theorem \ref{main}.
Let $\pi\colon X\to S$ be a flat Calabi-Yau fibration of
relative dimension two,
take a polarisation $\ell$ on $X$
and let $\hpi\colon Y\to S$ be a fine, irreducible component of
$\M^{\ell}(X/S)$.
Assume $Y$ and $X$ have the same dimension, $n$
say, and let $\P$ be a universal sheaf on $Y\times X$.

The main problem is that we do not know {\it a priori} that $Y$ is
non-singular, so we cannot immediately apply Theorem \ref{borlov}.
We also do not have
any control of the $\Ext$ groups between sheaves supported on the
singular fibres of $X$. These problems can be solved by using
Theorem \ref{biggy}.\footnote{Note added in 2019. There is a small logical gap in the published version at this point. Since $Y$ is defined to be an irreducible, rather than connected, component of the moduli space, the corresponding family $\{\P_y:y\in Y\}$ need not be complete. However, examining the proof of Theorem \ref{biggy} we see that in fact all one needs is that the Kodaira-Spencer map at each point $y\in Y$ is injective, which is certainly the case. Thanks to Nick Addington and Dan Bragg for pointing this out.} 
To apply it, let $U$ be an open
subset of $S$ over which $\pi$ is smooth.
Let $B$ be the complement of $U$ in $S$, a closed
subset of $S$ of positive codimension.

Note that $\P_y\tensor\omega_X=\P_y$ for any point $y\in\hpi^{-1}(U)$
because $\omega_X$ restricts to give the trivial bundle on each non-singular
fibre of $\pi$. By semi-continuity, for any $y\in Y$, there is a non-zero map
\[\P_y\lRa{\theta}\P_y\tensor\omega_X.\]
It is a consequence of Definition \ref{def} that the restriction of
$\omega_X$ to any fibre of $\pi$ is a numerically trivial line bundle.
It follows that the sheaves $\P_y$ and $\P_y\tensor\omega_X$ are
stable with the same Hilbert polynomials, so the map $\theta$ is an
isomorphism.

Let $y_1$, $y_2$ be distinct points of $Y$. The groups
$\Ext^i_X(\P_{y_1},\P_{y_2})$
vanish unless $y_1$ and $y_2$ lie on the same fibre of $\hpi$
and the
argument above shows that these groups also vanish when
$y_1$ and $y_2$ both lie over a point of $U$.
Thus the subscheme $\Gamma(\P)$ of the theorem is contained in the
union of the diagonal in $Y\times Y$ with the closed subscheme
\[Y\times_S Y\times_S B=\{(y_1, y_2)\in Y\times Y:\hpi(y_1)=\hpi(y_2)\in
B\}.\]
It only remains to show that this scheme has dimension at most
$n+1$.
But this is clear since the dimension of $B$ is at most $n-3$
and we assumed that $\hpi$ is equidimensional so
the fibres of $Y\times_S Y\to S$ all have dimension 4.

In the case when
$\pi$ is an elliptic fibration, $B$ has
dimension at most $n-2$ and the fibres of $Y\times_S Y\to S$ have
dimension $2$, so the same argument works.

To show that the fibration $\hpi\colon Y\to S$ satisfies the
condition of Definition \ref{def}, take a curve $C\subset Y$,
contained in a fibre of $\hpi$, and put $E=\OO_C$.
Then the object $\Phi(E)$ is supported on a fibre of
$\pi$. It follows that $\Phi(E)\tensor\omega_X$ has the same
numerical invariants as $\Phi(E)$, so using Serre duality
\begin{eqnarray*}
\eu(E\tensor\omega_Y)&=&\eu(E,\OO_Y)=\eu(\Phi(E),\Phi(\OO_Y)) \\
&=&\eu(\Phi(\OO_Y),\Phi(E))=\eu(\OO_Y,E)=\eu(E).
\end{eqnarray*}
This implies that $\K_Y\cdot C=0$. The canonical bundle of the general
fibre $Y_s$ of $\hpi$ is trivial because $\Phi$ restricts to give a
Fourier-Mukai transform $\Phi_s\colon \D(Y_s)\lra\D(X_s)$.
This completes the proof of Theorem
\ref{main}.


\section{FM transforms for elliptic threefolds}

In this section we study FM transforms
for elliptic threefolds.
In many places we follow the treatment
of FM transforms for elliptic surfaces given in \cite{Br1},
but in higher dimensions several new ideas are required.
In particular we shall prove that Mukai dual fibrations of
flat (equidimensional) fibrations are themselves flat.

\subsection{}
Let us fix the following notation.

\begin{defn}
An elliptic threefold is
a Calabi-Yau fibration $\pi\colon X\to S$, as in Definition \ref{def},
such that $S$ and $X$ have dimensions two and three respectively.
\end{defn}

Let $\pi\colon X\to S$ be an elliptic threefold. Write
$f\in H^4(X,\Z)$ for the Chern
character of the structure sheaf of a non-singular fibre of $\pi$ and
$\tau\in H^6(X,\Z)$
for the Chern character of a skyscraper sheaf.
Given an object $E$ of $\D(X)$
we put
\[d(E)=\cl_1(E)\cdot f\]
and refer to this number as the \emph{fibre degree} of $E$.
Note that $d(E)$ is the degree of the restriction of $E$ to a general fibre
of $\pi$.
Also define $\la_{X/S}$ to be the
highest common factor of the fibre degrees of objects of $\D(X)$.
Equivalently $\la_{X/S}$ is the smallest positive integer such that there
is a divisor $\sigma$ on $X$ with $\sigma\cdot f=\la_{X/S}$.

\subsection{}
We shall need the following piece of birational geometry.

\begin{prop}
\label{birat}
Let $\pi_i\colon X_i\to S$ be two elliptic threefolds.
Suppose there is a birational equivalence
$X_1\dashrightarrow X_2$ commuting
with the maps $\pi_i$. Then $\pi_2$ is flat if $\pi_1$ is.
\end{prop}

\begin{pf}
Let $\pi\colon X\to S$ be an elliptic threefold. Let
$\Delta\subset S$ be
the discriminant locus of $\pi$. Thus $\pi$ is smooth over
the open subset $U=S\setminus\Delta$. By Bertini's theorem,
we can find a non-singular curve $C\subset S$ meeting each irreducible
component of $\Delta$ such
that $\pi^{-1}(C)$ is non-singular. Since the morphism
$\pi\colon\pi^{-1}(C)\to C$ is a relatively minimal elliptic surface,
there is a positive
integer $d$ such that the natural map
\[\eta\colon \pi^*\pi_* (\omega_X^{\tensor d})\lra \omega_X^{\tensor d}\]
is an isomorphism when restricted to $\pi^{-1}(C)\subset X$.
This map is also an isomorphism on $\pi^{-1}(U)$ and hence
is an isomorphism away
from a finite number of fibres of $\pi$.
It follows that there is a
$\QQ$-divisor $\Lambda$ on $S$ such that $\K_X=\pi^* \Lambda + D$ where
$\pi(D)$ has codimension 2. Thus $D$ is contained in the union of
those fibres of
$\pi$ which have dimension 2.

Write $K_{X_i}=\pi_i^*(\Lambda_i)+D_i$. If $\pi_1$ is flat
then $D_1=0$. A simple argument \cite[Lemma 1.5]{Gr} then shows that
some positive multiple of $D_2$ is an effective divisor.
If $C$ is a general
hyperplane section of this divisor then $\K_{X_2}\cdot C<0$,
which is impossible, so $D_2=0$. It is then easy to show that
$\pi_2$ is flat \cite[Theorem 2.4]{Gr}.
This completes the proof.
\qed
\end{pf}

\subsection{}
\label{outahere}
The following result is basically a restatement of Theorem \ref{main}.

\begin{thm}
\label{duff}
Let $\pi\colon X\to S$ be an elliptic threefold.
Take a polarization
$\ell$ on $X$ and
let $Y\subset\M^{\ell}(X/S)$ be a fine, irreducible component of dimension 3
which contains a sheaf supported on a non-singular fibre of $\pi$.
Let $\P$ be a
universal sheaf on $Y\times X$. Then $Y$ with
the induced morphism $\hpi\colon Y\to S$
is an elliptic threefold and the functor
$\Phi=\Phi^{\P}_{Y\to X}$
is a FM transform.
\end{thm}

\begin{pf}
The only point to note is that we do not require $\hpi$ to be equidimensional.
Since $Y$ is irreducible $\hpi$ can only have finitely many fibres of
dimension 2, so $Y\times_S Y$ has dimension 4 and the argument of
Section \ref{bl}
still applies.
\qed
\end{pf}

Let $\pi\colon X\to S$ be an elliptic threefold
and take a divisor $\sigma$ on $X$ with $\sigma\cdot f=\la_{X/S}$.
Fix a pair of integers $a$ and $b$ such that $a$ is positive and
$b$ is coprime to $a\la_{X/S}$.

If $b=0$ then $a=\la_{X/S}=1$. If $\Phi\colon\D(Y)\to\D(X)$
is a FM transform taking points of $Y$
to sheaves of Chern character $f-\tau$,
then composing with a twist by $\OO_X(\sigma)$ gives a FM transform
taking points of $Y$ to sheaves of Chern character $f$.

Thus we may suppose that $b$ is non-zero.
Replacing a polarisation $\ell$ of $X$ by $\sigma\pm bm\ell$
for $m\gg 0$ we may assume
that $\ell\cdot f$ is coprime to $b$.
Let $Y$ be the unique
irreducible component of $\M^{\ell}(X/S)$ containing a point which represents
a stable bundle of rank $a$ and degree $b$ supported on a non-singular
fibre of $\pi$. As we noted in Section \ref{sixone}, the space $Y$
has dimension
three and is fine, so Theorem \ref{duff} shows that there is a FM
transform $\Phi\colon\D(Y)\to\D(X)$ taking structure sheaves of
points of $Y$ to sheaves of Chern character $af+b\tau$ on $X$.

\subsection{}
Take notation and assumptions as in Theorem \ref{duff}.
Let $\Psi$ denote the functor
$\Phi^{\Q}_{X\to Y}$
where
\[\Q=\R\lHom_{\OO_{Y\times X}}(\P,\pi_X^*\omega_X)[n-1].\]
The inverse of $\Phi$ is the functor $\Psi[1]$.
The reason for the strange choice of shift in the definition of $\Q$
is the following
 
\begin{lemma}
The object $\Q$ is a sheaf on $Y\times X$, flat over $Y$.
If $\hpi$ is flat
then $\P$ and $\Q$ are also flat over $X$.
\end{lemma}

\begin{pf}
The sheaf $\P$ is flat over $Y$, so for every point $(y,x)\in Y\times X$
\[\Ext^i_{Y\times X}(\P,\OO_{(y,x)})=\Ext^i_X(\P_y,\OO_x).\]
But $\Hom_X(\OO_x,\P_y)=0$ because $\P_y$ has pure dimension 1, so by Serre
duality these groups vanish
unless $0\leq i\leq 2$. Thus the homological dimension of
$\P$ and hence also of $\Q$ is 2. But $\Q$ is supported on $Y\times_S X$
which has codimension 2 in $Y\times X$, so $\Q$ is a sheaf.

Given a point $y\in Y$ let $i_y\colon \{y\}\times X\embed Y\times X$ denote the
inclusion. The object $\Q_y=\L i_y^*(\Q)$ has homological dimension 2 and
its support (which is the same as that of $\P_y$) has dimension 1.
Thus $\Q_y$
is a sheaf and $\Q$ is flat over $Y$.
If $\hpi$ is flat, the same argument shows that $\P$ and $\Q$
are flat
over $X$.
\qed
\end{pf}

\begin{cor}
\label{al}
If $\hpi$ is flat then for any sheaf $E$ on $X$
and ample line bundle $L$ the sheaf $E\tensor L^n$ is
$\Psi$-WIT$_0$ for $n\gg 0$.
\end{cor}

\begin{pf}
By the last lemma $\Q$ is a sheaf on $Y\times X$, flat over $X$, so
\[\Psi^i(E\tensor L^n)=\R\pi_{Y,*}(\Q\tensor\pi_X^*(E\tensor
L^n)).\]
Since $\pi_X^* L$ is $\pi_Y$-ample these groups vanish for $i>0$ when
$n\gg 0$.
\qed
\end{pf}

The next lemma is proved by a simple base-change.

\begin{lemma}
\label{base-change}
Given a point $s\in S$ such that the fibre $X_s$
is non-singular, let $i_s\colon X_s\embed X$ and $j_s\colon Y_s\embed Y$
be the corresponding embeddings of schemes. Then
there is an isomorphism of functors
\[\L i_s^*\circ\Phi\isom\Phi_s\circ \L j_s^*\]
where $\Phi_s$ is the functor $\Phi^{\P_s}_{Y_s\to X_s}$ and $\P_s$,
the restriction of the sheaf $\P$ to $Y_s\times X_s$, is a universal
sheaf parameterising stable sheaves on the fibre $X_s$.
\qed
\end{lemma}

\subsection{}
We can now prove the main result of this section.

\begin{prop}
\label{equidimensional}
In the situation of Theorem \ref{duff}, the fibration
$\hpi\colon Y\to S$ is flat precisely if
$\pi\colon X\to S$
is.
\end{prop}

\begin{pf}
Suppose first that $\hpi$ is flat. Fix an
ample line bundle $L$ on $X$. Each sheaf $\P_y$ is
supported in dimension 1, so for all $y\in Y$ and any integer $i\neq 1$,
$\Ext^i_X(L^n,\P_y)=0$ for $n\gg 0$.
Since
\[\Hom_{\D(Y)}^i(\Psi(L^n),\OO_y)=\Hom_{\D(X)}^{i+1}(L^n,\P_y)\]
it follows that for $n\gg 0$ the object $\Psi(L^n)$ is a
locally-free sheaf on $Y$.

Given a point $s\in S$
we can use Corollary \ref{al}
to obtain a $\Psi$-WIT$_0$ sheaf $F$ on $X$
whose support is the fibre
$X_s$ of $\pi$. But 
\[\Ext^2_X(L^n,F)=\Ext^2_Y(\Psi(L^n),\Psi(F))=0 \mbox{ for }n\gg 0\]
because $\Psi(F)$ is supported on the fibre $Y_s$ of $\hpi$ which has
dimension 1. This implies that the support of $F$ has dimension 1 and so
$\pi$ is equidimensional, hence flat.

Now suppose that $\pi$ is flat. For any point $s\in S$ such that
the fibre $X_s$ is non-singular the restriction of $\P$ to $Y_s\times X_s$
is a universal sheaf parameterising stable bundles on $X_s$. It follows
from the results of \cite[Section 3]{Br1}
that for any $x\in X$ lying on a non-singular fibre of $\pi$
the sheaf $\P_x$ is a stable bundle on
the elliptic curve $Y_s$ of rank $a$
and degree $c$
coprime to $a\la_{Y/S}$.
As in Section \ref{outahere}
we can construct a new elliptic threefold
$\pi'\colon Z\to S$ parameterising sheaves
of Chern character $af+c\tau$ on $Y$.
The fibrations $\pi\colon X\to S$ and $\pi'\colon Z\to S$ are then birational.
Proposition \ref{birat}
shows that $\pi'$ is flat, so what we proved above shows that
$\hpi\colon Y\to S$ is also flat. This completes the proof.
\qed
\end{pf}


\section{Further properties of the transforms}

In this section we show that the properties of FM transforms
for elliptic surfaces listed in \cite[Section 6]{Br1} also hold in the case
of flat elliptic threefolds. As an application we prove Theorem
\ref{moduli}.

\subsection{}
Take notation as in the last section.
Putting together what we proved there gives
the following generalisation of \cite[Theorem 5.3]{Br1}.

\begin{thm}
\label{super}
Let $\pi\colon X\to S$ be a flat elliptic threefold.
Take an element
\[\mat{c}{a}{d}{b}\in\SL_2(\Z),\]
such that $a>0$ and $\la_{X/S}$ divides $d$. Then there is a
a flat elliptic threefold
$\hpi\colon Y\to S,$
and a sheaf $\P$ on $Y\times X$,
supported on $Y\times_S X$ and flat over both factors,
such that for any point
$(y,x)\in Y\times X$, $\P_y$ has Chern character $af+b\tau$ on $X$ and
$\P_x$ has Chern character $af+c\tau$ on $Y$.

The resulting functor 
$\Phi=\Phi^{\P}_{Y\to
X}$ is a FM transform satisfying
\begin{equation}
\label{num}
\col{r(\Phi E)}{d(\Phi E)}=\mat{c}{a}{d}{b}\col{r(E)}{d(E)},
\end{equation}
for all objects $E$ of $\D(Y)$.
\end{thm}

\begin{pf}
The fibration $\hpi\colon Y\to S$ is constructed as in Section
\ref{outahere}.
Lemma \ref{base-change} and the results of \cite[Section 2]{Br1}
show that there exist integers
$c$ and $d$ such that (\ref{num}) holds. It follows that $\la_{X/S}=\la_{Y/S}$.
Twisting $\P$ by pull-backs of  line bundles from $Y$
we obtain all possible values of
$c$ and $d$.
\qed
\end{pf}

\subsection{}
Take notation and assumptions as in Theorem \ref{super}.
Since $\P$ is flat over $Y$,
with support of  dimension one over $X$, for any sheaf $E$ on $Y$
one has
\[\Phi^i(E)=0 \mbox{  unless }0\leq i\leq 1.\]
Furthermore $\Phi$ is left-exact so that for any short exact sequence
\[0\lra A\lra B\lra C\lra 0\]
on $X$ one obtains a long exact sequence
\begin{eqnarray*}
0 &\lra & \Phi^0(A) \lra \Phi^0(B) \lra \Phi^0(C) \\
  &\lra & \Phi^1(A) \lra \Phi^1(B) \lra \Phi^1(C) \lra 0.
\end{eqnarray*}
Similar statements hold for $\Psi$.

\smallskip

The isomorphism
$\Psi\circ\Phi=\id_{\D(Y)}[-1]$ implies that if $E$ is a $\Phi$-WIT$_i$ sheaf on
$X$ then $\Hat{E}$ is a $\Psi$-WIT$_{1-i}$ sheaf on $Y$. More generally,
for any sheaf $E$ on $Y$ there is
a spectral sequence
\[ E^{p,q}_2=\Psi^p(\Phi^q(E)) \implies \left\{\begin{array}{ll}
E &\mbox{ if $p+q=1$} \\
0 &\mbox{ otherwise.}
\end{array} \right.\]
>From what we showed above one has $E^{p,q}_2=0$ unless $p=0$ or
$1$ so the spectral sequence degenerates. This leads to the following
results.

\smallskip

\begin{lemma}
\label{seqq}
Let $E$ be a sheaf on $Y$. Then $\Phi^0(E)$ is $\Psi$-WIT$_1$, $\Phi^1(E)$
is $\Psi$-WIT$_0$, and there is a short exact sequence
\[0\lra \Psi^1(\Phi^0(E))\lra E\lra \Psi^0(\Phi^1(E))\lra 0. \qed\]
\end{lemma}

\smallskip

\begin{lemma}
\label{wit}
A sheaf $E$ on $Y$ is $\Phi$-WIT$_0$ if and only if
$\Hom_Y(E,\Q_x)=0$
for all $x\in X$.
\end{lemma}

\begin{pf}
See \cite[Lemma 6.5]{Br1}.
\qed
\end{pf}

\begin{lemma}
\label{tf}
If a torsion-free sheaf $E$ on $Y$ is $\Phi$-WIT$_0$, the transformed sheaf
$\Hat{E}$
is also torsion-free.
\end{lemma}

\begin{pf}
We repeat the argument of \cite[Lemma 7.2]{Br1}.
Note first that the sheaf $\Hat{E}$ is $\Psi$-WIT$_1$.
If $T$ is a torsion subsheaf we obtain a short exact sequence
\[0\lra T\lra \Hat{E}\lra Q\lra 0\]
and applying $\Psi$ gives an exact sequence
\[0\lra\Psi^0(Q)\lra\Psi^1(T)\lRa{g} E\lra\Psi^1(Q)\lra0.\]
Thus $T$ is $\Psi$-WIT$_1$. Since $T$ is a torsion sheaf, $\cl_1(E)$
is the Chern class of an effective divisor, so
for any point $y\in Y$
\[\eu(\Hat{T},\OO_y)=\eu(T,\P_y)=-\cl_1(T)\cdot af\leq 0.\]
It follows that $\Hat{T}$ has rank 0. Since $E$ is torsion-free
the map $g$ is zero, so $\Psi^0(Q)=\Hat{T}$. But the first
sheaf is $\Phi$-WIT$_1$ while the second is $\Phi$-WIT$_0$,
so both sheaves are zero and hence $T=0$.
\qed
\end{pf}

\begin{lemma}
\label{stab}
If $E$ is a $\Phi$-WIT$_0$ sheaf on $Y$ whose restriction
to the general fibre of
$\hpi$ is stable then the restriction of the transform $\Hat{E}$
to the general fibre of $\pi$ is stable.
\end{lemma}

\begin{pf}
Lemma \ref{base-change} and \cite[Proposition 3.3]{Br1}.
\qed
\end{pf}

\subsection{}
We can now prove Theorem \ref{moduli} which
is restated below.

\medskip

\paragraph*{{\bf Theorem \ref{moduli}.}\ }
{\it Let $X$ be a threefold with a flat elliptic fibration
$\pi\colon X\to S$ and let $\hpi\colon Y\to S$ be a Mukai dual fibration
as in Theorem \ref{main}. Let $\NN$ be a connected component of
the moduli space of rank one, torsion-free sheaves on $Y$.
Then there is a polarisation $\ell$ of $X$ and a connected component
$\M$ of the moduli space of stable torsion-free sheaves on $X$
with respect to $\ell$
which is isomorphic to $\NN$.}

\medskip

\begin{pf}
Given a point $s\in \NN$, let $\E_s$ be the corresponding rank 1,
torsion-free on $Y$. The argument of Lemma \ref{al} shows that
twisting by a line bundle
on $Y$ we may
assume that all these sheaves are $\Phi$-WIT$_0$.
By Lemmas \ref{tf}, \ref{stab} and \ref{fib} we can find a
polarisation
$\ell$ of $X$ such that each of the transformed sheaves
$\Hat{\E_s}$ is stable
with respect to $\ell$. Applying Proposition \ref{jo}
completes the proof.
\qed
\end{pf}

In fact we can say slightly more.

\begin{lemma}
\label{bundle}
Take assumptions as in Theorem
\ref{moduli} and suppose also that $\NN$ parameterises
ideal sheaves of subschemes $C\subset Y$
which meet each fibre of $\hpi$ in
at most a finite number of points.
Then under the isomorphism of the theorem,
locally-free sheaves in $\M$ correspond precisely to
equidimensional Cohen-Macaulay
subschemes $C\subset Y$ of dimension one.
\end{lemma}

\begin{pf}
For each point $s\in \NN$, the corresponding torsion-free sheaf
$\E_s$ on $Y$ is of the form $M\tensor\I_C$, where $M$
is a line bundle and  $C\subset Y$ is a subscheme of $Y$.
As in the proof of Theorem \ref{moduli} we can assume that $M$
is chosen so that the sheaves $\E_s$ are all $\Phi$-WIT$_0$.
There is an isomorphism
\[\Ext^i_X(\Hat{\E_s},\OO_x)=\Ext^{i+1}_Y(\E_s,\Q_x)\]
and $\Hat{\E_s}$ is locally-free if and only if these groups vanish for all
$i>0$ and all $x\in
X$. By Serre duality
\[\Ext^3_Y(\E_s,\Q_x)=\Hom_Y(\Q_x,\E_s\tensor\omega_Y)^{\dual}=0\]
because $\E_s$ is torsion-free. Applying the functor
$\Hom_Y(-,\Q_x)$ to the exact sequence
\[ 0\lra \E_s\lra M\lra M\tensor\OO_C\lra 0\]
and noting that $\Ext^i_Y(M,\Q_x)=0$ for $i>1$ shows that
\[\Ext^2_Y(\E_s,\Q_x)=\Ext^3_Y(M\tensor\OO_C,\Q_x).
=\Hom_Y(\Q_x,M\tensor\OO_C)^{\dual}\]
Since the supports of $\Q_x$ and $\OO_C$ intersect in a finite number
of points, this group is zero unless there is a non-zero map
$\OO_y\to\OO_C$ for some $y\in Y$.
\qed
\end{pf}

Taking $\NN$ to be a component of the moduli space of rank one, torsion-free
sheaves which parameterise
ideal sheaves of zero-dimensional subschemes of
$Y$ shows that there are connected components of the moduli space
of stable, torsion-free sheaves on $X$ which contain no locally-free sheaves.


\section{An example}

Let $X\subset \PP^2\times\PP^2$ be a non-singular $(3,3)$ (anticanonical)
divisor. We shall suppose that $X$ does not contain a fibre of either
of the projections $p_i\colon \PP^2\times\PP^2\to \PP^2$. This is equivalent to
the assumption that each of the induced projections $\pi_i\colon X\to\PP^2$
is flat. Thus $X$ is a Calabi-Yau threefold with two
flat elliptic fibrations $\pi_1$, $\pi_2$.
In this section we apply the results of the previous section to
compute a moduli space of stable bundles on $X$.

\subsection{}
By the Lefschetz hyperplane theorem
$H^1(X,\C)=0$ and $H^2(X,\C)=\C^{\oplus 2}$ is generated by the two
elements $H_i=\pi_i^*(L)$, where $L$ is the cohomology class
of a line on $\PP^2$. By
Poincar{\'e} duality, $H^4(X,\C)=\C^{\oplus 2}$ is generated by the
two elements $f_i=\pi_i^*(\omega)$, where $\omega$ is the cohomology class
of a point on $\PP^2$. We have relations
\[ H_1^2=f_1, \quad H_1 \cdot H_2=f_1+f_2, \quad H_2^2=f_2,
 \quad H_1\cdot f_2= H_2\cdot f_1=3. \]

We shall consider $X$ to be elliptically fibred via the map $\pi_1$. Thus
in the notation of the previous section we take $\pi=\pi_1$ and $f=f_1$.
One can check directly that $\pi_*(\OO_X)$ and $\R^1\pi_*(\OO_X)$ are
locally-free
and it follows that for any fibre $X_s$ of $\pi$ one has
\[H^0(X,\OO_{X_s})=H^1(X,\OO_{X_s})=\C.\]

Let $\Delta\subset X\times X$ denote the diagonal
and let $\P$ be the kernel of the restriction map
$\OO_{X\times_S X}\to\OO_{\Delta}$. Thus there is a short exact
sequence
\begin{equation}
\label{bj}
0\lra \P\lra \OO_{X\times_S X}\lra \OO_{\Delta}\lra 0.
\end{equation}
Considered as a sheaf on $X\times X$, $\P$ is flat over both
factors, and for any point $x\in X$ there is a short exact sequence
\begin{equation}
\label{sp}
0\lra \P_x\lra \OO_{X_s}\lra \OO_x\lra 0,
\end{equation}
where $X_s$ is the scheme-theoretic fibre of $\pi$ containing the
point $x$. Put $\Phi=\Phi^{\P}_{X\to X}$.

\begin{lemma}
\label{ally}
For any object $E$ of $\D(X)$ there is a triangle in $\D(X)$
\[ \Phi(E)\lra \pi^* \R\pi_*(E)\lRa{\eta} E\lra\Phi(E)[1],\]
where $\eta$ is the unit of the adjunction $\pi^*\vdash\R\pi_*$.
\end{lemma}

\begin{pf}
The short exact sequence
(\ref{bj}) implies that there is a triangle in $\D(X)$
\[ \Phi(E)\lra \Theta(E)\lra E\lra \Phi(E)[1]. \]
Here $\Theta$ is the functor
\[ \Theta(-)=\Phi^{\OO_{X\times_S X}}_{X\to X}=\R q_{2,*}\circ q_1^* \]
where $q_1$, $q_2$ are the projections $X\times_S X\to X$. Base-changing
around the diagram
\[
\begin{CD}
X\times_S X &@> q_1 >> &X \\
@V q_2 VV   &          &@V \pi VV \\
X           &@>\pi>>   &S
\end{CD}
\]
shows that in fact $\Theta(-)=\pi^* \R\pi_*(-)$.
\qed
\end{pf}

As before, let $\Psi[1]$ be a left adjoint for $\Phi$.
Then $\Psi\isom\Phi^{\Q}_{X\to X}$ where $\Q$ is a
sheaf on $X\times X$, flat over both factors. For each $x\in X$ there is
a short exact sequence
\[0\lra \OO_{X_s}\lra \Q_x\lra \OO_x\lra 0.\]

\begin{lemma}
The functor $\Phi$ is a Fourier-Mukai transform.
\end{lemma}

\begin{pf}
It is not clear that the sheaf $\P_x$ is stable for all $x\in X$ so we
shall show that $\Phi$ is a FM transform directly.
Since $\Psi(\OO_x)=\Q_x$ it will be enough to check that for all $x\in
X$, $\Phi(\Q_x)=\OO_x[-1]$.

By the short exact sequence (\ref{sp}), $H^0(X,\P_x)=0$ and $H^1(X,\P_x)=\C$.
Then Serre duality gives $H^0(X,\Q_x)=\C$ and $H^1(X,\Q_x)=0$, so
$\R\pi_*(\Q_x)=\OO_s$ where $s=\pi(x)$.
Thus there is a triangle
\[\Phi(\Q_x)\lra \OO_{X_s}\lRa{\eta} \Q_x\lra \Phi(\Q_x)[1]\]
where the map $\eta$ is non-zero and hence injective.
The result follows.
\qed
\end{pf}

\subsection{}
Let $\NN$ denote the connected component of the moduli space
of rank one, torsion-free sheaves on $X$
containing the ideal sheaf of a fibre of $\pi_2$.
There is an obvious morphism
$\beta\colon\PP_2\to\NN$ taking
a point in $\PP^2$ to the ideal sheaf of the fibre of $\pi_2$ lying over it.

\begin{lemma}
The morphism $\beta$ is an
isomorphism.
\end{lemma}

\begin{pf}
The morphism $\beta$ sends a point $t\in T=\PP^2$ to the ideal
sheaf $E_t=\pi_2^*(\I_t)$. It is clearly injective on points.
Using the fact that $\pi_2$ is flat it is easy to check that
$\Ext^1_X(E_t,E_t)$ is two-dimensional
for all $t\in T$, so $\beta$ is an isomorphism.  
\qed
\end{pf} 

Let $M$ be the line bundle $\OO_X(H_2)$.
Given an ideal sheaf $\I_C\in\NN$ put $E=M\tensor\I_C$. Note
that $M\tensor\OO_C=\OO_C$ so there is a short exact sequence 

\begin{equation}
\label{jjf}
0\lra E\lra M\lra \OO_C\lra 0.
\end{equation}

\begin{lemma}
The sheaf $E$ is $\Phi$-WIT$_0$.
\end{lemma}

\begin{pf}
By Lemma \ref{wit} we must show that $\Hom_X(E,\Q_x)=0$ for all $x\in X$.
If $\Q_x$ is supported on the fibre $D$ of $\pi$,
any map $E\to \Q_x$ factors
through $E|_D$.
There is an exact sequence
\[0\lra \Tor_1^{\OO_X}(\OO_C,\OO_D)\lra
E|_D\lRa{g} M\tensor\OO_D\lra \OO_C\tensor\OO_D\lra
0.\]
Since $\Q_x$ has no zero-dimensional subsheaves, any map $E|_D\to\Q_x$ factors
through the image of $g$.

Note that $\OO_C\tensor\OO_D$ is either
the structure sheaf of a point or zero, depending on whether $C$ and $D$
meet or not.
Thus it will be enough to check that for any $y\in X$ there
are no non-zero maps
$M\tensor\P_y\to\Q_x.$

Let $A$ be the image of such a map.
The fact that $A$ is a subsheaf of $\Q_x$ implies that $H^0(X,A)$ has
dimension at most 1 and  so $\eu(A)\leq 1$.
Similarly $H^1(X,A\tensor M^{-1})$ has dimension at most 1, and
since $M$ has fibre degree 3 this gives
$\eu(A)\geq 2$, a contradiction.
\qed
\end{pf}

\subsection{}

We now compute the Chern character of the transformed sheaf $\Hat{E}$. To do this,
note first that Lemma \ref{ally} implies that
\[\ch(\Hat{E})=\pi^*\ch(\R\pi_*(E))-\ch(E).\]
The sequence (\ref{jjf}) together with the fact that $\eu(\OO_C)=0$ give
\[ \ch(E)=(1,H_2,\frac{1}{2}f_2,0)-(0,0,f_2,0)=(1,H_2,-\frac{1}{2} f_2,0). \]

Let $i\colon X\embed \PP^2\times\PP^2$ be the embedding.
Applying the functor $\R p_{1,*}$
to the short exact sequence of sheaves on $\PP^2\times\PP^2$
\[0\lra \OO_{\PP^2\times\PP^2}(-3,-2)\lra
\OO_{\PP^2\times\PP^2}(0,1)\lra i_* M\lra 0\]
shows that $\R \pi_{1,*}(M)=\OO_S^{\oplus 3}$.

Considered as a subscheme of $\PP^2\times\PP^2$, $C$ is of the form
$\{x\}\times D$ with $D\subset \PP^2$ a curve of degree 3 and
$\R p_{1,*}(\OO_C)=\OO_D$.
Thus $\R\pi_{1,*}(E)$ has Chern character
\[ (3,0,0)-(0,3,-\frac{9}{2})=(3,-3,\frac{9}{2}) \]
and $\Hat{E}$ has Chern character
\[ (3,-3H_1,\frac{9}{2} f_1,0)-
(1,H_2,-\frac{1}{2} f_2,0)=(2,-3H_1-H_2,\frac{9}{2} f_1 +\frac{1}{2} f_2,0). \]
To tidy this up, let $L$ be the line bundle $\OO_X(2H_1+H_2)$. Then
\[\ch(L)=(1,2H_1+H_2,4f_1+\frac{5}{2}f_2,9)\]
so that $\Hat{E}\tensor L$ has Chern character
$(2,H_1+H_2,\frac{3}{2}f_1-\ha f_2,0).$

By Lemmas \ref{stab} and \ref{bundle}
the sheaf $\Hat{E}$ is locally-free with
stable restriction to the general fibre of $\pi$.
Applying Lemmas \ref{fib}, \ref{hilb} and Proposition \ref{jo} gives

\begin{prop}
There is a polarisation $\ell$ on $X$ and a connected component
of the moduli space of stable sheaves with respect to $\ell$
of Chern character 
\[(2,H_1+H_2,\frac{3}{2}f_1-\ha f_2,0)\]
which is isomorphic to $\PP^2$.
Moreover, all elements of this component are $\mu$-stable vector
bundles.
\qed
\end{prop}


\bigskip

\noindent Department of Mathematics and Statistics,
The University of Edinburgh,
King's Buildings, Mayfield Road, Edinburgh, EH9 3JZ, UK.

\smallskip

\noindent email: {\tt tab@maths.ed.ac.uk}
{\tt  \ \ \  A.Maciocia@ed.ac.uk}

\end{document}